\documentclass[a4paper,10pt]{article}
\def \Z {{\mathbf {Z}}}

\usepackage[T1]{fontenc} % РёРЅРѕРіРґР° T1 РІРјРµСЃС‚Рѕ T2A
\usepackage[cp1251]{inputenc}
\usepackage[russian]{babel}

\begin{document}
\Large
\vspace{20mm}
{\bf Asymmetries in  Asymptotic  3-fold Properties of Ergodic Actions}
\\
\vspace{20mm}

We present: 
\vspace{10mm}

1)  a mixing $Z ^ 2$-action  with   the following asymmetry of multiple mixing property:
for some commuting measure-preserving   transformations $S$, $T$ and a sequence $n_j$ 
$$
\lim_{j\to  \infty}\mu(A\bigcap S^{-n_j}A\bigcap T^{-n_j}A)=\mu(A)^3
  $$
 for all measurable sets $A$, but there is   $A_0$,  $\mu(A_0)=\frac 1 2$, such that      
$$
\lim_{j\to  \infty}\mu(A_0\bigcap S^{n_j}A_0\bigcap T^{n_j}A_0)=0;
  $$  
\vspace{10mm}

2)   $Z $-actions with the asymmetry of the partial multiple  mixing and the partial multiple rigidity:
$$
\lim_{j\to \infty}\mu(A\bigcap T^{k_j}A\bigcap T^{m_j}A)=
\frac23 \mu(A)^3+\frac13\mu(A),
 $$
$$
\lim_{j\to \infty}\mu(A\bigcap T^{-k_j}A\bigcap T^{-m_j}A)=
 \mu(A)^2;
 $$
\vspace{10mm}

3) an infinite transformation $T$ such that  for all $A$, $\mu(A)<\infty$,
$$
\lim_{j\to \infty}\mu(A\bigcap T^{k_j}A\bigcap T^{m_j}A)=
\frac13\mu(A)
 $$
and
$$
\lim_{j\to \infty}\mu(A\bigcap T^{-k_j}A\bigcap T^{-m_j}A)=0.
  $$
\newpage

\begin{center}
{   \bf Об асимметрии кратных асимптотических свойств \\  эргодических  действий
}
\end{center}
 
\vskip .1in

\begin{center}{Рыжиков В.В.}
\end{center}

\vskip .1in

{\large Аннотация.  В статье предъявлены перемешивающие $\Z^2$-действия, не изоморфные своим обратным; рассмотрены $\Z$-действия  с асимметрией  свойств  частичного кратное перемешивания на последовательностях и частичной кратной жесткости;  приведены новые примеры преобразований  пространства с бесконечной мерой, не изоморфных своим обратным.  }

\Large
\vskip .1in
\section{ Асимметрия  кратного перемешивания для  $\Z^2$-действий.} В работе ~\cite{1}  В.А. Рохлин ввел понятие кратного перемешивания, предполагая, что новый метрический инвариант  сможет различать динамические системы с одинаковым спектром.  
До сих пор не доказано существование перемешивающих преобразований, не обладающих кратным перемешиванием. Однако, для $\Z^2$-действий  Ф. Ледраппье  
 ~\cite{2} дал изящное решение  $\Z^2$-аналога этой знаменитой  проблемы.   Оказалось, что коммутирующие автоморфизмы некоторых компактных коммутативных групп  относительно меры  Хаара дают    примеры перемешивающих действий, не обладающих кратным перемешиванием. 
Действия типа  Ледраппье, как мы покажем, могут обладать  асимметрией свойства кратного перемешивания.

\textbf{Теорема 1.} \it Существует  перемешивающее  действие, порожденное коммутирующими   преобразованиями $S$, $T$, такое, что для некоторой последовательности $n_m\to +\infty$  для всех измеримых множеств  $A$  
$$
\lim_{m\to  \infty}\mu(A\bigcap S^{-n_m}A\bigcap T^{-n_m}A)=\mu(A)^3,
  $$
но для  некоторого   множества $A_0$ меры $\frac 1 2$  для всех $n_m$ выполнено 
$$
\mu(A_0\bigcap S^{n_m}A_0\bigcap T^{n_m}A_0)=0.
  $$
Иначе говоря, найдется $\Z^2$-действие, которое  кратно перемешивает относительно  пары последовательностей $(-n_m,0), \ (0,-n_m)$  при $m\to\infty$, а обратное действие
 этим  свойством не обладает. 
\rm

Доказательство.   В качестве фазового пространства действия выбирается
группа  $X$    всех последовательностей $a:\Z^2\to\Z_2$, удовлетворяющих тождеству
$$a(z_1,z_2)+a(z_1+1,z_2)+a(z_1,z_2+1) =0.$$
Групповая операция -- почленное  сложение (mod 2) последовательностей.
 Важную роль для дальнейшего играет  наблюдение Ледраппье о том, что исходное тождество влечет за собой серию других тождеств:
$$a(z_1,z_2)+a(z_1+2^m,z_2)+a(z_1,z_2+2^m) =0. \eqno (\ast)$$
 Действительно,  из
$$a(z_1,z_2)+a(z_1+1,z_2)+a(z_1,z_2+1) +a(z_1+1,z_2)+a(z_1+2,z_2)+a(z_1+1,z_2+1)+$$  
$$+a(z_1,z_2+1)+a(z_1+1,z_2+1)+a(z_1,z_2+2)=a(z_1,z_2)+a(z_1+2,z_2)+a(z_1,z_2+2)$$
имеем
$$0+0+0= a(z_1,z_2)+a(z_1+2,z_2)+a(z_1,z_2+2).$$
Теперь аналогично  получаем  $a(z_1,z_2)+a(z_1+2,z_2)+a(z_1,z_2+4)=0$  и т.д.

Рассмотрим автоморфизмы $$ Ta(z_1,z_2)=a(z_1+1,z_2), \ \  Sa(z_1,z_2)=a(z_1,z_2+1).$$
Мера Хаара $\mu$  инвариантна относительно  $S$ и $T$.
 Хорошо известно, что действие, порожденное этими автоморфизми обладает  перемешиванием.

Докажем, что для любых нетривиальных характеров $\chi_1$, $\chi_2$, $\chi_3$  для всех достаточно больших $m$  выполнено
$$\int_X\chi_1 (S^{-2^m}\chi_2) (T^{-2^m} \chi_3)d\mu =0, \eqno (0)
$$
что непосредственно  влечет за собой
$$\mu(A\cap S^{-2^m}B\cap T^{-2^m}C)\to \mu(A)\mu(B)\mu(C), \ \ m\to+\infty.
$$
так как характеры образуют базис в $L_2(\mu)$.

Пусть $(0)$  не выполняется, тогда для бесконечного множества
чисел $m$
$$\chi_1 (S^{-2^m}\chi_2) (T^{-2^m} \chi_3) =1.
$$
В силу $(\ast)$   выполняется
$$a+S^{2^m}a+T^{2^m}a=0,$$ $$  \chi (S^{2^m}\chi)  (T^{2^m} \chi)=1 \eqno (1)$$
для любого характера $\chi$. Подставляя $\chi=S^{-2^m}T^{-2^m} \chi_3$ в $(1)$,      
получим
$$S^{-2^m}T^{-2^m} \chi_3 (T^{-2^m} \chi_3)  (S^{-2^m} \chi_3)=1. $$
Умножение  левой  части этого выражения  на  $\chi_1 (S^{-2^m}\chi_2) (T^{-2^m} \chi_3)=1$  дает
$$
\chi_1 \chi_3(S^{-2^m}\chi_2\chi_3) (S^{-2^m}T^{-2^m} \chi_3) =1.
$$
Теперь умножим  на 
$$
S^{-2^m}(\chi_1 \chi_3) \ (\chi_1\chi_3)\  S^{-2^m}T^{2^m} (\chi_1\chi_3) =1
$$
и получим
$$
S^{-2^m}(\chi_1 \chi_2) \   S^{-2^m}T^{2^m} (\chi_1\chi_3)(S^{-2^m}T^{-2^m} \chi_3)=1,
$$
что равносильно
$$
\chi_1 \chi_2 \   T^{2^m} (\chi_1\chi_3)T^{-2^m} \chi_3 =1.
$$
Как  доказал Рохлин \cite{1},\S 3,  эргодический  автоморфизм $T$  компактной коммутативной группы
перемешивает кратно. Следовательно, последнее тождество
возможно для больших значений $m$  только 
в случае 
$$
\chi_1 \chi_2 =\chi_1\chi_3= \chi_3 =1,
$$
следовательно, $
\chi_1 =\chi_2 =\chi_3 =1.
$
Итак, первая часть утверждения теоремы доказана.

Для доказательства второй части  мы  предъявим множество  $A_0$,   $\mu(A_0)=\frac 1 2$, такое, что 
$
\mu(A_0\bigcap S^{2^m}A_0\bigcap T^{2^m}A_0)=0
  $
   для всех  $m>0$.
Для этого фиксируем нетривиальный характер $\chi$ и обозначим через
 $A_0$  прообраз $\{ -1\}$ относительно функции $\chi$.

 Заметим, что  характеры вида $\chi  S^{2^m}\chi $, \ $S^{2^m}\chi  T^{2^m}\chi$ нетривиальны (это вытекает из эргодичности автоморфизмов $ S^{2^m}$ и $S^{-2^m}T^{2^m}$).  Пользуясь $(1)$, получим  для всех  $m>0$
$$\int_X(1- \chi )(1-S^{2^m}\chi ) (1- T^{2^m} \chi )d\mu$$ 
$$ =1- 
\int_X \chi \ d\mu +\dots +  \int_X  S^{2^m}\chi \ T^{2^m} \chi d\mu -
\int_X \chi \ S^{2^m}\chi \ T^{2^m} \chi d\mu$$  
$$=1 -0 +\dots +0 -1=0.$$ 
Так как подынтегральная функция 
$(1- \chi )(1-S^{2^m}\chi ) (1- T^{2^m} \chi )$ неотрицательна, получим, что  она равна 0 тождественно. Поэтому 
$
\mu(A_0\bigcap S^{2^m}A_0\bigcap T^{2^m}A_0)=0.
  $

Отметим, что рассмотренное действие и его обратное обладают 
счетнократным лебеговским спектром, см. \cite{3}.

\section{ Асимметрия частичного кратного перемешивания и 
частичной кратной жесткости преобразований} В  ~\cite{4}, ~\cite{5} для обнаружения эффекта асимметрии $\Z$-действий  использовались  
такие свойства как частичная кратная жесткость и частичное кратное перемешивание на  последовательностях.  Эти инварианты
(свойства, инвариантные  относительно сопряжения)  позволяли   различать действие и его обратное. Недавно И.С. Ярославцев    упростил  метод работы  ~\cite{2}. Процитируем основное утверждение  заметки ~\cite{6}.

 \it Существует эргодическое преобразование $T$, обладающее  свойством: для некоторой последовательности $n_j\to\infty$ и для любых  измеримых множеств $A, B$ и $C$ выполнено равенство
$$
\lim_{j\to \infty}\mu(A\bigcap T^{n_j+1}B\bigcap T^{2n_j}C)=
$$
 $$ =\frac13 \left(\mu(A\bigcap TB\bigcap T^{-1}C)+\mu(A\bigcap T^{2}B\bigcap TC)+\mu(A\bigcap B\bigcap C) \right).
 $$
\rm

При наличии такого предела преобразования $T$ и $T^{-1}$ не изоморфны друг другу\cite{4}.
Это вытекает из того, что   предел 
$\lim_{j\to \infty}\mu(A\bigcap T^{-n_j-1}B\bigcap T^{-2n_j}C)$  не содержит компоненты $\mu(A\bigcap B\bigcap C)$. Более подробные пояснения см. ниже.
\rm

Повторяя рассуждения из  ~\cite{4},  получаем следующее формально более общее утверждение.

\textbf{Теорема 2.} \it Существует эргодическое преобразование $T$, обладающее следующим свойством: для любого $N>0$ найдется   последовательность $n_j\to\infty$ такая, что  для любых  измеримых множеств $A, B$ и $C$ выполнено 

$$
\lim_{j\to \infty}\mu(A\bigcap T^{n_j+N}B\bigcap T^{2n_j}C)=
$$

 $$ =\frac13 \left(\mu(A\bigcap T^NB\bigcap T^{-N}C)+\mu(A\bigcap T^{2N}B\bigcap T^NC)+\mu(A\bigcap B\bigcap C) \right).
 $$
\rm
 
 \bf Замечание. \rm Из приведенного свойства преобразования  $T$ вытекает, что оно обладает слабым перемешиванием.  Действительно, оператор 
$\frac13 \left(T^{-1}+ T+ I \right)$ является слабым пределом некоторой последовательности степеней преобразования, поэтому $T$ не имеет собственных функций, кроме констант. Действительно,
если $|\lambda|=1$ и $\lambda^{n_i}\to \frac13 (\lambda^{-1} +\lambda +1)$,  то  $\lambda=1$.  Для эргодического $T$ только постоянные функции удовлетворяют условию $Tf=f$.  Отсутствие нетривиальных
собственных функций означает слабое перемешивание (непрерывность спектра).

{\bf Конструкция преобразования.} Примеры преобразований, удовлетворяющих условиям теоремы 2,   строятся в классе преобразований ранга 1 (их определение см., например, в  \cite{7}).  Достаточно потребовать, чтобы для каждого
$N>0$  нашлась бесконечная последовательность этапов построения преобразования ранга 1, в каждом из  которых
массив надстроек (spacers) имел бы вид
$(0,2N,N,0,2N,N,...,0,2N,N)$, причем длины этих массивов не ограничены. 
Тогда по аналогии с  
работой  \cite{6} для   такой конструкции $T$ при любом фиксированном $N$ найдутся последовательность
  $h_j\to\infty$
и  последовательность разбиений  фазового пространства
$$
  X= Y_i\sqcup Y^1_i \sqcup Y^2_i\sqcup  Y^3_i
$$
такие, что  $\mu(Y_i)\to 0$,  а  множества
$Y^1_i,Y^2_i,Y^3_i$ имеют следующий вид:
\medskip
\medskip

$
  Y^1_j = B^1_j\cup TB^1_j\cup  \dots
\cup T^{h_j-2}B^1_j\cup T^{h_j-1}B^1_j,
$
\medskip

$
  Y^2_j = B^2_j\cup TB^2_j\cup \dots
  \cup T^{h_j}B^2_j  \cup T^{h_j+2N -1}B^2_j,
$
\medskip

$
  Y^3_j = B^3_j  \cup TB^3_j \cup \dots
   \cup T^{h_j-1}B^3_j  \cup T^{h_j+N-1}B^3_j.
$
\\
Также потребуем выполнения  условий: 
 последовательность разбиений  $$
 \{ C_j, TC^1_j, T^2C_j, \dots , T^{h_j-1}C_j\},\
  \     
C_j= B^1_j\cup B^2_j\cup B^3_j,
$$
стремится к разбиению на точки,
$$
  T^{h_j}B^1_j= B^2_j, \ \ 
 T^{h_j+2N}B^2_j= B^3_j,
\ \ 
\mu(T^{h_j+N}B^3_j\Delta B^1_j)/\mu(B^1_j) \ \to \ 0.
$$
 Эти условия  обеспечивают 
наличие  пределов, фигурирующих в теореме 2. Отметим, что  в заметке \cite{4} был рассмотрен случай $N=1$.

\section{ Доказательство асимметрии } Предварительно установим равенство 
$$
\lim_{j\to \infty}\mu(A\bigcap T^{-n_j-N}B\bigcap T^{-2n_j}C)=$$
$$
= \frac13 \left(\mu(A\bigcap B\bigcap T^{N}C)+\mu(A\bigcap T^{-N}B\bigcap T^{-N}C)+\mu(A\bigcap T^{-2N}B\bigcap C) \right).$$
Заметим, что 
$$
\lim_{j\to \infty}\mu(A\bigcap T^{-n_j-N}B\bigcap T^{-2n_j}C)=
\lim_{j\to \infty}\mu(T^{2n_j}A\bigcap T^{2n_j}T^{-n_j-N}B\bigcap T^{-2n_j}T^{-2n_j}C)
$$
$$=\lim_{j\to \infty}\mu(C\bigcap T^{n_j-N}B\bigcap T^{2n_j}A)=
\lim_{j\to \infty}\mu(C\bigcap T^{n_j+N}T^{-2N}B\bigcap T^{2n_j}A)
$$
и применим к последнему пределу теорему 2, подставив $T^{-2N}B$  вместо $B$.  
Получаем, что   предел  $\lim_{j\to \infty}\mu(A\bigcap T^{-n_j-N}B\bigcap T^{-2n_j}C)$  равен  
$$ \frac13 \left(\mu(C\bigcap T^{N}T^{-2N}B\bigcap T^{-N}A)+\mu(C\bigcap B\bigcap T^NA)+\mu(C\bigcap T^{-2N}B\bigcap A) \right),
$$
что нам и   нужно.

Так как $T$ --  слабо перемешивающее преобразование,
 найдется последовательность $N_j\to\infty$  такая, что 
$$ \lim_j \mu(A\bigcap T^{N_j}B\bigcap T^{-N_j}C)=
 \lim_j \mu(A\bigcap T^{2N_j}B\bigcap T^{N_j}C)=\mu(A)\mu(B)\mu(C).
 $$
Это вытекает из известной   теоремы Фюрстенберга ~\cite{8},  утверждающей
 для большинства $n$ близость величин
$  \mu(A\bigcap T^{n}B\bigcap T^{2n}C)$ и  
 $\mu(A)\mu(B)\mu(C).$  
(Отметим, что можно не использовать упомянутую теорему, а построить явную конструкцию
с подходящей перемешивающей последовательностью $N_j$.) 
Пользуясь   теоремой 2,  в выражении
$ \mu(A\bigcap T^NB\bigcap T^{-N}C)+\dots+\mu(A\bigcap B\bigcap C) 
 $ переходим к (повторному)  пределу по $N$ вдоль упомянутой последовательности $N_j$.  
Получим
 $$ \lim_j \mu(A\bigcap T^{N_j}B\bigcap T^{-N_j}C)+\mu(A\bigcap T^{2N_j}B\bigcap T^{N_j}C)+\mu(A\bigcap B\bigcap C)\ $$
 $$=2\mu(A)\mu(B)\mu(C) + \mu(A\bigcap B\bigcap C).$$
 Таким образом,  
преобразование $T$  обладает частичным кратным перемешиванием 
относительно пары некоторых  последовательностей  (соответствующий предел имеет компоненту $\mu(A)\mu(B)\mu(C)$).  Но  преобразование $T^{-1}$ этим свойством не обладает. Действительно, имеем
$$\lim_j \mu(A\bigcap B\bigcap T^{N_j}C)+\mu(A\bigcap T^{-N_j}B\bigcap T^{-N_j}C)+\mu(A\bigcap T^{-2N_j}B\bigcap C) \ $$
$$
=\mu(A\bigcap B)\mu(C) +\mu(B\bigcap C)\mu(A) +\mu(A\bigcap C)\mu(B). 
$$  
Положив  $k_j = n_j +N_j$, $m_j=2n_j$ 
(рост перемешивающей  последовательности  $N_j$ выбирается достаточно  медленным), получим следующую асимметрию:
$$  \lim_{j\to \infty}\mu(A\bigcap T^{k_j}B\bigcap T^{m_j}C)=
\frac23 \mu(A)\mu(B)\mu(C)+\frac13\mu(A\bigcap B\bigcap C),$$
$$
\lim_{j\to \infty}\mu(A\bigcap T^{-k_j}B\bigcap T^{-m_j}C)=
\frac13 \left(\mu(A\bigcap B)\mu(C)+\mu(B\bigcap C)\mu(A)+\mu(A\bigcap C)\mu(B)\right).
 $$
Подведем итог в случае $A=B=C$.
%\newpage

\textbf{Теорема 3.} \it Существует слабо перемешивающее  преобразование $T$ такое, что для некоторых   последовательностей $k_j,m_j\to\infty$   и любого  измеримого множества $A$ выполнено 

$$
\lim_{j\to \infty}\mu(A\bigcap T^{k_j}A\bigcap T^{m_j}A)=
\frac23 \mu(A)^3+\frac13\mu(A),
 $$

$$
\lim_{j\to \infty}\mu(A\bigcap T^{-k_j}A\bigcap T^{-m_j}A)=
 \mu(A)^2.
 $$
\rm

\textbf{Следствие.} \it Преобразования $T$ и $T^{-1}$ не изоморфны.
\rm

Действительно,  асимптотические свойства в теореме 2 
являются инвариантами, т. е. сохраняются при сопряжении.
Если $S=R^{-1}TR$ и $R$ сохраняет меру, то
$$
\mu(A\bigcap S^{-k_j}A\bigcap S^{-m_j}A)=
\mu(RA\bigcap T^{-k_j}RA\bigcap T^{-m_j}RA).
$$
Поэтому
$$ \lim_{j\to \infty}\mu(A\bigcap S^{-k_j}A\bigcap S^{-m_j}A)=\mu(RA)^2=\mu(A)^2.
 $$
Аналогично устанавливается  равенство
$$ \lim_{j\to \infty}\mu(A\bigcap S^{k_j}A\bigcap S^{m_j}A)=\frac23 \mu(A)^3+\frac13\mu(A).
 $$

Таким образом,    показано,  что  \it пара последовательностей $k_j, m_j$ обеспечивает  для 
некоторого преобразования  одновременно и
частичное кратное перемешивание  (слагаемое $\frac23 \mu(A)^3$)  и частичную кратную жесткость  (слагаемое $\frac13 \mu(A)$),  в то время как  обратное преобразование 
этими свойствами не обладает.  \rm

\section{ Случай бесконечной меры.  Замечания}  Асимметрия сохраняющих меру преобразований особенно наглядна в случае, когда фазовое пространство имеет бесконечную
(сигма-конечную) меру.  Отметим, что свойство перемешивания, т. е.  условие $\mu(A\bigcap T^{N_j}B)\to\mu(A)\mu(B)$,  теперь  превращается
в $\mu(A\bigcap T^{N_j}B)\to 0$  для любых множеств $A,B$ конечной меры.    Приведем "бесконечный" аналог теоремы 3.
\vspace{3mm}

\textbf{Теорема 4.} \it Существует эргодическое  обратимое  преобразование $T$,
сохраняющее бесконечную меру,  такое, что для некоторых   последовательностей $k_j,m_j\to\infty$   и любого  измеримого множества $A$ конечной меры выполнено 

$$
\lim_{j\to \infty}\mu(A\bigcap T^{k_j}A\bigcap T^{m_j}A)=
\frac13\mu(A),
 $$

$$
\lim_{j\to \infty}\mu(A\bigcap T^{-k_j}A\bigcap T^{-m_j}A)=0.
  $$
\rm
\vspace{2mm}

 Доказательство этой теоремы осуществляется по тому же плану, что и  доказательства теоремы 3. Но теперь ситуация  упрощается за счет того, что перемешивание в пространстве с бесконечной мерой тривиальным  образом влечет за собой кратное перемешивание.    Конструкции   преобразований отличаются от упомянутых ранее конструкций на пространстве с конечной мерой лишь тем, что вместо последовательности массивов надстроек 
вида  $$(0,2N,N,0,2N,N,...,0,2N,N)$$  используются массивы  вида
$$(H_j,H_j+2N,H_j+N,H_j,H_j+2N,H_j+N,...,H_j,H_j+2N,H_j+N).$$ 
Последовательность $H_j$ выбирается так, чтобы мера фазового пространства конструкции оказалась бесконечной.  Здесь подойдет любая достаточно быстро растущая последовательность. 

В заключение отметим, что приведенные выше результаты о $\Z$-действиях можно перенести на потоки.
Поток, не изоморфный обратному,  действующий на
 пространстве бесконечной меры был предъявлен  в \cite{9}.
Изучению асимметрии асимптотических свойств потоков  посвящена  статья \cite{10}.  Наша работа, как и заметка \cite{6},
была стимулирована  поиском наиболее простых и наглядных инвариантов, которые показывают отличие действия и его обратного.  Для эргодической теории  представляет интерес отсутствия сплетений между различными положительными степенями преобразования. Из  результатов  работ по этой тематике (см., например, \cite{9}, \cite{11}, \cite{12}) видно, что  
неизоморфизм положительных степеней  может быть следствием  некоторых двукратных асимптотических свойств. Так как преобразования   $T$ и $T^{-1}$ обладают одинаковыми асимптотическими свойствами кратности 2,   необходимо
использовать свойства большей кратности, чтобы установить их неизоморфизм.
\large

\vspace{5mm}

Московский государственный  университет 

им. М.В. Ломоносова

vryzh@mail.ru

07.02.2014
\end{document}